\documentclass{amsart}

\newtheorem{theorem}{\textbf{Theorem}}
\newtheorem{corollary}{\textbf{Corollary}}

\newtheorem{conjecture}{\textbf{Conjecture}}
\newtheorem{lemma}{\textbf{Lemma}}
\newtheorem{exe}{\textbf{Example}}

\def\Z {\mathbb{Z}}
\def\Q {\mathbb{Q}}

\theoremstyle{remark}

\numberwithin{equation}{section}

\begin{document}

\title[SOME $U_m$-NUMBERS RELATED TO LIOUVILLE CONSTANT]{TRANSCENDENCE MEASURES FOR SOME $U_m$-NUMBERS RELATED TO LIOUVILLE CONSTANT}
\author{ANA PAULA CHAVES}
\address{DEPARTAMENTO DE MATEM\'{A}TICA, UNIVERSIDADE FEDERAL DO CEAR\' A, FORTALEZA, CEAR\' A, BRAZIL}
\email{apchaves.math@gmail.com}

\author{DIEGO MARQUES}
\address{DEPARTAMENTO DE MATEM\'{A}TICA, UNIVERSIDADE DE BRAS\'ILIA, BRAS\'ILIA, DF, BRAZIL}
\email{diego@mat.unb.br}




\begin{abstract}
In this note, we shall prove that the sum and the product of an algebraic number $\alpha$ by the \textit{Liouville constant} $L=\sum_{j=1}^{\infty}10^{-j!}$ is a $U$-number with type equals to the degree of $\alpha$ (with respect to $\mathbb{Q}$). Moreover, we shall have that 
\begin{center}
$\max\{w^{\ast}_n(\alpha L),w^{\ast}_n(\alpha + L)\}\leq 2m^2n+m-1$, for $n=1,...,m-1$.
\end{center}
\end{abstract}

\maketitle

\section{Introduction}
A real number $\xi$ is called a \textit{Liouville number}, if for any positive real number $w$ there exist infinitely many rational numbers $p/q$, with $q\geq 1$, such that
\begin{center}
$0<\left|\xi - \dfrac{p}{q}\right|<\dfrac{1}{q^w}$.
\end{center}

Transcendental number theory began in 1844 when Liouville \cite{lio} showed that all Liouville numbers are transcendental establishing thus the first examples of such numbers. For instance, the number
\begin{center}
$L=\displaystyle\sum_{j=1}^{\infty}10^{-j!}=0.11000100000000000000001000...$,
\end{center}
which is known as \textit{Liouville's constant}, is a Liouville number and therefore transcendental. In 1962, Erd$\ddot{\mbox{o}}$s \cite{erdos} proved that every nonzero real number can be written as the sum and the product of two Liouville numbers. 

In 1932, Mahler \cite{mah} splited the set of the transcendental numbers in three disjoint sets named \textit{$S$-,\ $T$-} and {\it $U$-numbers}. Particularly, the $U$-numbers generalizes the concept of Liouville numbers. We denote by $w^*_n(\xi)$ as the supremum of the real numbers $w^*$ for which there exist infinitely many real algebraic numbers $\alpha$ of degree $n$ satisfying
\begin{center}
$0<|\xi - \alpha| < H(\alpha)^{-w^*-1}$,
\end{center}
where $H(\alpha)$ (so-called the \emph{height} of $\alpha$) is the maximum of absolute value of coefficients of the minimal polynomial of $\alpha$ (over $\Z$). The number $\xi$ is said to be a \emph{$U^*_m$-number} (according to LeVeque \cite{le}) if $w^*_m(\xi)=\infty$ and $w^*_n(\xi)<\infty$ for $1\leq n< m$ ($m$ is called the \emph{type} of the $U$-number). We point out that we actually have defined a Koksma $U^*_m$-number instead of a Mahler $U_m$-number. However, it is well-known that they are the same \cite[cf. Theorem 3.6]{bugeaud} and \cite{al2}. We observe that the set of $U_1$-numbers is precisely the set of Liouville numbers.

The existence of $U_m$-numbers, for all $m\geq 1$, was proved by LeVeque \cite{le}. In 1993, Pollington \cite{pol} showed that for any positive integer $m$, every real number can be expressed as a sum of two $U_m$-numbers.

Since two algebraically dependent numbers must belong to the same Mahler's class \cite[Theorem 3.2]{bugeaud}, then $\alpha L$ and $\alpha+L$ are $U$-numbers, for any nonzero algebraic number $\alpha$. But what are their types?

In this note, we use G$\ddot{\mbox{u}}$tting's method \cite{gut} for proving that the sum and the product of every $m$-degree algebraic number $\alpha$ by $L$ is a $U_m$-number. Moreover, we obtain an upper bound for $w_n^{\ast}$.

\begin{theorem}
Let $\alpha$ be an algebraic number of degree $m$ and let $L$ be the Liouville's constant. Then $\alpha L$ and $\alpha+L$ are $U_m$-numbers, with
\begin{equation}\label{in}
\max\{w^{\ast}_n(\alpha L),w^{\ast}_n(\alpha + L)\}\leq 2m^2n+m-1,\ for\ n=1,...,m-1.
\end{equation}
\end{theorem}

\section{Auxiliary Results}
Before the proof of the main result, we need two technical results. The first one follows as an easy consequence of the triangular inequality and binomial identities.
\begin{lemma}\label{1}
Given $P(x)\in \Z[x]$ with degree $m$ and $a/b\in \Q \backslash \{0\}$. If $Q_1(x)=a^mP(\frac{b}{a}x)$ and $Q_2(x)=b^mP(x-\frac{a}{b})$, then
\begin{itemize}
\item[\textup{(i)}] $H(Q_1)\leq \max\{|a|,|b|\}^mH(P)$;
\item[\textup{(ii)}] $H(Q_2)\leq 2^{m+1}\max\{|a|,|b|\}^mH(P)$.
\end{itemize}
\end{lemma}
\begin{proof}
(i) If $P(x)=\sum_{j=0}^ma_jx^j$, then $Q_1(x)=\sum_{j=0}^ma_jb^ja^{m-j}x^j$. Supposing, without loss of generality, that $|a|\geq |b|$, we have $|a|^m|a_j|\geq |a|^{m-j}|a_j||b|^j$ for $0\leq j\leq m$. Hence, we are done. For (ii), write $Q_2(x)=\sum_{i=0}^mc_ix^i$, where 
\begin{center}
$c_i=b^m\displaystyle\sum_{j=i}^ma_j{j\choose j-i}(-1)^{j-i}\left(\frac{a}{b}\right)^{j-i}$
\end{center}
Therefore
\begin{center}
$|c_i|\leq H(P)\displaystyle\sum_{k=0}^{m-i}\displaystyle{k+i\choose k}|a|^k|b|^{m-k}\leq \max\{|a|,|b|\}^mH(P)\displaystyle\sum_{k=0}^{m-i}\displaystyle{k+i\choose k}$.
\end{center}
Since $\sum_{k=0}^{m-i}{k+i\choose k}={m+1\choose m-i}\leq 2^{m+1}$, we finally have
\begin{center}
$|c_i|\leq 2^{m+1}\max\{|a|,|b|\}^mH(P)$,
\end{center}
which completes our proof.
\end{proof}

In addition to Lemma 1, we use the fact that algebraic numbers are not well aproximable by algebraic numbers.
\begin{lemma}[Cf. Corollary A.2 of \cite{bugeaud}]\label{2}
Let $\alpha$ and $\beta$ be two distinct nonzero algebraic numbers of degree $n$ and $m$, respectively. Then we have
\begin{center}
$|\alpha-\beta|\geq (n+1)^{-m/2}(m+1)^{-n/2}\max\{2^{-n}(n+1)^{-(m-1)/2},2^{-m}(m+1)^{-(n-1)/2}\}$\\
$\times H(\alpha)^{-m}H(\beta)^{-n}$.
\end{center}
\end{lemma}
\begin{proof}
A sketch of the proof can be found in the Appendix A of \cite{bugeaud}.
\end{proof}

\section{Proof of the Theorem}

For $k\geq 1$, set
\begin{center}
$p_k=10^{k!}\displaystyle\sum_{j=1}^k10^{-j!}$, $q_k=10^{k!}$ and $\alpha_k=\dfrac{p_k}{q_k}$.
\end{center}
We observe that $H(\alpha_{k-1})<H(\alpha_k)=10^{k!}=H(\alpha_{k-1})^{k}$ and 
\begin{equation}\label{lp}
|L - \alpha_k| < \frac{10}{9}H(\alpha_k)^{-k-1}.
\end{equation}
Thus, setting $\gamma_k=\alpha \alpha_k$, we obtain of (\ref{lp})
\begin{equation}\label{11}
|\alpha L - \alpha \alpha_k|\leq cH(\alpha_k)^{-k-1},
\end{equation}
where $c=10|\alpha|/9$. It follows by the Lemma \ref{1} (i) that $H(\alpha_k)^m\geq H(\alpha)^{-1}H(\gamma_k)$ and thus we conclude that
\begin{equation}\label{opo}
|\alpha L - \alpha \alpha_k|\leq cH(\alpha)^{k+1}H(\gamma_k)^{-k-1}.
\end{equation}
Consequently, $\alpha \beta$ is a $U$-number with type at most $m$ (since $\gamma_k$ has degree $m$).

Again, we use Lemma \ref{1} (i) for obtaining
\begin{equation}
H(\gamma_{k+1})\leq H(\alpha)H(\alpha_{k+1})^m=H(\alpha)H(\alpha_k)^{(k+1)m}\leq H(\alpha)H(\gamma_k)^{(k+1)m}
\end{equation}

Now, let $\gamma$ be an $n$-degree real algebraic number, with $n<m$ and $H(\gamma)\geq H(\gamma_1)$. Thus, one may ensure the existence of a sufficient large $k$ such that
\begin{equation}\label{key}
H(\gamma_k)< H(\gamma)^{2m^2}< H(\gamma_{k+1})\leq H(\alpha)H(\gamma_k)^{(k+1)m}.
\end{equation}
So, by Lemma \ref{2}, it follows that
\begin{equation}\label{eq1}
|\gamma_k - \gamma|\geq f(m,n)H(\gamma)^{-m}H(\gamma_k)^{-n},
\end{equation}
where $f(m,n)$ is a positive number which does not depend on $k$ and $\gamma$. Therefore by (\ref{key}) 
\begin{equation}\label{22}
|\gamma_k - \gamma|\geq f(m,n)H(\alpha)^{-1/2m}H(\gamma_k)^{-(k+1)/2-n}.
\end{equation}
By taking $H(\gamma)$ large enough, the index $k$ satisfies
\begin{equation}\label{33}
H(\gamma_k)^{(k+1)/2-n}\geq 2cf(m,n)^{-1}H(\alpha)^{k+1/2m}.
\end{equation}
Thus, it follows from (\ref{opo}), (\ref{22}) and (\ref{33}) that $|\gamma_k-\gamma|\geq 2|\alpha L - \gamma_k|$. Therefore, except for finitely many algebraic numbers $\gamma$ of degree $n$ strictly less than $m$, we have
\begin{eqnarray*}
|\alpha L - \gamma| & \geq & |\gamma_k - \gamma|-|\alpha L-\gamma_k|\nonumber\\
 & \geq & \frac{1}{2}|\gamma_k - \gamma| \nonumber\\
 & \geq & \frac{f(m,n)}{2}H(\gamma)^{-m}H(\gamma_k)^{-n} > \dfrac{f(m,n)}{2}H(\gamma)^{-2m^2n-m},\nonumber\\
\end{eqnarray*}
where we use the left-hand side of (\ref{key}). It follows that $w_n^*(\alpha L)\leq 2m^2n+m-1$ which finishes our proof.

The case $\alpha+L$ follows the same outline, where we use Lemma \ref{1} (ii) rather than (i).
\qed

\section{The general case and further comments}

Let $\beta$ be a Liouville number. Since that a $U$-number keeps its type when multiplied by any nonzero rational number, we can consider $0<\beta < 1$. Set
\begin{center}
$S_{\beta}=\{(\frac{p_k}{q_k})_{k\geq 1}\in \mathbb{Q}^{\infty}\ :\ |\beta-\frac{p_k}{q_k}|<\frac{1}{q_k^{k+1}},\ k=1,2,...\}$.
\end{center}
By the assumption on $\beta$, we may suppose $1\leq p_k\leq q_k$ and then $H(p_k/q_k)=q_k$, for all $k$. Note that $S_{\beta}$ is an infinite set. 

As is customary, the symbols $\ll$, $\gg$ mean that there is an implied constant in the inequalities $\leq$, $\geq$, respectively. In our process for proving the Theorem 1, the key step happens when holds an inequality like in (\ref{key}). Thus it follows that
\begin{theorem}
Let $\alpha$ be an $m$-degree algebraic number and let $\beta$ be a Liouville number. If there exists a sequence $(p_k/q_k)_{k\geq 1} \in S_{\beta}$ such that $q_k\ll q_{k+1}\ll q_{k}^{k+1}$ for all $k\gg 1$, then the numbers $\alpha \beta$ and $\alpha+\beta$ are $U_m$-numbers and
\begin{equation}\label{in2}
\max\{w^{\ast}_n(\alpha \beta),w^{\ast}_n(\alpha + \beta)\}\leq 2m^2n+m-1,\ for\ n=1,...,m-1.
\end{equation}
\end{theorem}
\begin{exe}
For any integer number $m\geq 2$ and any $a_j\in \{1,...,9\}$, the number $\sum_{j=1}^{\infty}a_jm^{-j!}$ is a Liouville number satisfying the hypothesis of the previous theorem.
\end{exe}
\begin{corollary}
For any $m\geq 1$, there exists an uncountable collection of Liouville numbers that are expressible as sum of two algebraically dependent $U_m$-numbers.
\end{corollary}
\begin{proof}
Set $\beta=\displaystyle\sum_{j=1}^{\infty}a_j10^{-j!}$, where $a_j\in \{1,2\}$. The result follows immediately of Theorem 2 and of writing $\beta = (\frac{\beta+\sqrt[m]{2}}{2})+(\frac{\beta-\sqrt[m]{2}}{2})$.
\end{proof}

There exist several lower estimates for the distance between two distinct algebraic numbers, e.g., Liouville's inequality and Lemma \ref{2}. A too-good-to-be-true Conjecture due to Schmidt \cite{sch} states that
\begin{conjecture}
For any number field $\mathbb{K}$ and any positive real number $\epsilon$, we have
\begin{center}
$|\alpha-\beta| > c(\mathbb{K},\epsilon)(\max\{H(\alpha),H(\beta)\})^{-2-\epsilon}$,
\end{center}
for any distinct $\alpha,\beta \in \mathbb{K}$, where $c(\mathbb{K},\epsilon)$ is some constant depending only on $\mathbb{K}$ and on $\epsilon$.
\end{conjecture}

We conclude by pointing that if the Schmidt's conjecture is true, then the sum and the product of any $m$-degree algebraic number $\alpha$ by any Liouville number $\beta$ is a $U_m$-number and the inequality (\ref{in2}) can be considerable improved for 
\begin{center}
$\max\{w^{\ast}_n(\alpha \beta),w^{\ast}_n(\alpha + \beta)\}\leq 1$.
\end{center}


\section*{Acknowledgement}
The second author thanks to the Department of Mathematics of Universidade Federal do Cear\' a for its hospitality. We also thank Yann Bugeaud for nice discussions on the subject.




\end{document}